\documentclass[12pt]{article}
\usepackage{dsfont}
\usepackage{amsfonts}
\usepackage{amsthm}
\usepackage{amssymb}
\usepackage{mathtools}
\usepackage{amsmath}
\usepackage{mathrsfs}
\usepackage{mathrsfs}
\usepackage{stmaryrd}
\usepackage{framed}
\usepackage{nccmath}
\usepackage{wrapfig}
\usepackage{enumitem}
\usepackage{geometry}
\usepackage{tikz,tikz-cd}
\usepackage{bm}
\usepackage[stable]{footmisc}
\usepackage[new]{old-arrows}
\usepackage{draftwatermark}
\usepackage{quiver}
\SetWatermarkText{DRAFT - GSN - 2022}
\SetWatermarkColor[gray]{0.95}
\SetWatermarkScale{0.4}
\NewDocumentCommand{\tens}{t_}
 {%
  \IfBooleanTF{#1}
   {\tensop}
   {\otimes}%
 }
\NewDocumentCommand{\tensop}{m}
 {%
  \mathbin{\mathop{\otimes}\displaylimits_{#1}}%
 }

\newcommand{\prtt}[1]{\left( #1 \right)}
\newcommand{\tA}[2]{#1 \tens_{R} #2}

\newcommand{\bC}{{\mathbb{C}}}

\newcommand{\bN}{{\mathbb{N}}}

\newcommand{\sub}{\subseteq}
\newcommand{\ten}{\otimes}
\newcommand{\fp}{\mathfrak{p}}
\newcommand{\fq}{\mathfrak{q}}
\newcommand{\fm}{\mathfrak{m}}

\newcommand{\id}{\mbox{id}}

\newcommand{\spec}{\mbox{Spec}}

\theoremstyle{plain}

\newtheorem{defi}{Definition}[section]
\numberwithin{defi}{section}

\newtheorem{prop}[defi]{Proposition}
\newtheorem{teo}[defi]{Theorem}
\newtheorem{cor}[defi]{Corollary}
\newtheorem{lema}[defi]{Lemma}
\newtheorem{example}[defi]{Example}

\begin{document}
	
	\title{Universally injective and integral contractions on relative Lipschitz saturation of algebras}
	\author{Thiago da Silva and Maico Ribeiro}
	\date{}
	
	\maketitle
	
	\begin{abstract}
		In this work, we obtain contraction results for a class of diagrams of ring morphisms which strictly includes the ones obtained by Lipman. Further, we present some applications on quotient and in the changing of the base ring in the saturation.
	\end{abstract}
	
	\let\thefootnote\relax\footnote{2020 \textit{Mathematics Subjects Classification} 13B22, 14B07
		
		}

	\section*{Introduction}
	
	%The study in equisingularity was started by Zariski in \cite{zariski1965}, where he was interested to investigate this concept in an algebraic variety along an irreducible (singular) subvariety. First, he dealt with different ways to define equivalent singularities of plane algebroid curves, and after in \cite{zariski1965-2} Zariski worked with algebroid hypersurfaces with a singular point at its origin.
	
	%In the meantime, in \cite{zariski1968} introduced an operation on a ring $A$, which he called \textit{saturation}, which consists in passing from $A$ to some ring $\tilde{A}$ lying over $A$ and under the integral closure of $A$ in its total ring of fractions. In that work Zariski showed how useful this operation is for the theory of singularities by means of geometric applications to plane algebroid curves, and more generally, to algebroid hypersurfaces.
	
	%Having its applications in his mind, firstly Zariski restricted himself to the case in which $A$ is a local domain of dimension one, and always with the assumption that the base field is algebraically closed and of characteristic zero. Then, in \cite{zariski1971(0), zariski1971, zariski1975} he presented his general theory of saturation, extending several results and showing how to apply them to a more general setup in the theory of singularities.
	
	The core of saturation theory developed by Zariski \cite{zariski1965, zariski1965-2, zariski1968, zariski1971(0), zariski1971, zariski1975} is to look for a special intermediate algebra between a ring $A$ and its integral closure $\overline{A}$. In \cite{phamt}, in the case of complex analytic algebras, Pham and Teissier observed that the germs of Lipschitz meromorphic functions lie between $A$ and $\overline{A}$, showing that it coincides with the Zariski saturation in the hypersurface case. Namely, for a reduced complex analytic algebra $A$ with normalization $\overline{A}$, Pham and Teissier defined the Lipschitz saturation of $A$ as $$A^*:=\{f\in\overline{A}\mid f\ten_{\bC}1-1\ten_{\bC}f\in\overline{I_A}\},$$
	
	\noindent where $I_A$ denotes the kernel of the canonical map $\overline{A}\ten_{\bC}\overline{A}\rightarrow \overline{A}\ten_A\overline{A}$. Although Pham and Teissier were thinking with an analytic background, they left a good question about if $A^*_B\sub \overline{A}$, in the case where $B$ is an $A$-algebra and $$A^*_B:=\{f\in B\mid f\ten_{\bC}1-1\ten_{\bC}f\in\overline{I_{A, B}}\},$$
	
	\noindent where $I_{A,B}$ now denotes the kernel of the canonical map $B\ten_{\bC}B\rightarrow B\ten_AB$. 
	
	In \cite{lipman1} Lipman extended this definition for a sequence of ring morphisms $R\rightarrow A\overset{g}{\rightarrow}B$, and defined what he called the \textit{relative Lipschitz saturation of $A$ in $B$}, denoted by $A^*_{B, R}$. Besides,  Lipman developed several techniques and general properties on this operation in the ring $A$. One of them was to investigate if the sequence of ring morphisms % https://q.uiver.app/#q=WzAsNCxbMCwwLCJSIl0sWzEsMCwiQSJdLFsyLDAsIkIiXSxbMywwLCJCJyJdLFswLDFdLFsxLDIsImciXSxbMiwzLCJmIl1d
	\[\begin{tikzcd}
		R & A & B & {B'}
		\arrow[from=1-1, to=1-2]
		\arrow["g", from=1-2, to=1-3]
		\arrow["f", from=1-3, to=1-4]
	\end{tikzcd}\] 
	
	\noindent satisfies the pullback condition  $f^{-1}(A^*_{B',R})=A^*_{B,R}$, which was proved to be true provided $f$ is faithfully flat or $f$ is integral with $\ker f$ nil ideal.

	In this work, we provide more general conditions than those presented by Lipman in \cite{lipman1} to guarantee this pullback condition holds, working on diagrams that generalize these sequences. In this way, in Section 1 we recall some important definitions and we define what we call a \textit{Lipman diagram} (see Definition \ref{202306211340}), which is a special type of diagram of ring morphisms on the form 
	
	% https://q.uiver.app/#q=WzAsNixbMCwwLCJSIl0sWzEsMCwiQSJdLFsyLDAsIkIiXSxbMCwxLCJSJyJdLFsxLDEsIkEnIl0sWzIsMSwiQiciXSxbMCwxXSxbMSwyLCJnIl0sWzMsNF0sWzQsNSwiZyciXSxbMCwzXSxbMSw0XSxbMiw1XV0=
	\[\begin{tikzcd}
		R & A & B \\
		{R'} & {A'} & {B'}
		\arrow[from=1-1, to=1-2]
		\arrow["g", from=1-2, to=1-3]
		\arrow[from=2-1, to=2-2]
		\arrow["{g'}", from=2-2, to=2-3]
		\arrow[from=1-1, to=2-1]
		\arrow[from=1-2, to=2-2]
		\arrow[from=1-3, to=2-3]
	\end{tikzcd}\]

	In the main theorem of Section 2, we prove that Lipman diagrams contract (i.e, $f^{-1}((A')^*_{B', R'})=A^*_{B, R}$) if $f$ and $R\rightarrow R'$ are universally injective, which is a strictly more general condition than faithful flatness, which was the condition required in the related Lipman result in \cite{lipman1}.
	
	In Section 3, we show that Lipman diagrams contract only assuming that $f$ is integral, with no assumptions about $\ker f$ being a nil ideal, as required by Lipman in \cite{lipman1}.
	
	Further, at the end of sections 2 and 3, we apply Theorem \ref{202306011654} and Theorem \ref{202306021356} to get new classes of algebras where the changing of the base ring $R\rightarrow  R'$ does not affect the relative Lipschitz saturation of $A$ in $B$ in a sequence % https://q.uiver.app/#q=WzAsNCxbMCwwLCJSIl0sWzEsMCwiUiciXSxbMiwwLCJBIl0sWzMsMCwiQiJdLFswLDFdLFsxLDJdLFsyLDMsImciXV0=
	\[\begin{tikzcd}
		R & {R'} & A & B
		\arrow[from=1-1, to=1-2]
		\arrow[from=1-2, to=1-3]
		\arrow["g", from=1-3, to=1-4]
	\end{tikzcd}\]

	\section{Lipman diagrams}

	Let us remind some basic definitions and fix important notations. Consider a sequence of ring morphisms 
	
	\begin{align*}
		R \overset{\tau}{\longrightarrow} A \overset{g}{\longrightarrow} B.
	\end{align*}
\noindent which induces algebra structure on $A$ and $B$. Consider the map from $B$ to its tensor product with itself over $R$
	\begin{align*}
		\Delta : B &\longrightarrow \tA{B}{B} \\
		b &\longmapsto b\ten_R1_B-1_B\ten_Rb.
	\end{align*}
	
	The following properties of $\Delta$ is an easy exercise which we may safely leave to the reader. 
	
	\begin{enumerate}
		\item $\Delta(b_1+b_2)=\Delta(b_1+b_2), \forall b_1,b_2\in B$;
		
		\item $\Delta(rb)=r\Delta(b), \forall r\in R$ and $b\in B$;
		
		\item (Leibniz rule) $\Delta(b_1b_2)=(b_1\ten_R 1)\Delta(b_2)+(1\ten_Rb_2)\Delta(b_1), \forall b_1,b_2\in B$. 
	\end{enumerate}
	
	By the universal property of the tensor product, there exists a unique $R$-algebras morphism $$\varphi: \tA{B}{B} \longrightarrow B\tens_{A}B$$
	
	\noindent which maps $x\otimes_R y\mapsto x\otimes_A y$, for all $x,y\in B$, that is known as the \textit{canonical morphism}. Furthermore, we already know (see \cite{kleiman}, 8.7) that $$\ker\varphi=\langle ax\ten_R y-x\ten_Ray\mid a\in A\mbox{ and }x,y\in B \rangle.$$
	
	Notice that $ax\ten_R y-x\ten_Ray=(x\ten_Ry)(g(a)\ten_R1-1\ten_Rg(a))$. Therefore,$$\ker\varphi=\langle g(a)\ten_R1-1\ten_Rg(a)\mid a\in A\rangle=\Delta(g(A))(B\ten_RB),$$
	
	\noindent i.e, $\ker\varphi$ is the ideal of $B\ten_RB$ generated by the image of $\Delta\circ g$. Now we recall the definition of the relative Lipschitz saturation.
	
	\begin{defi}\cite{lipman1}
		The Lipschitz saturation of $A$ in $B$ relative to $R \overset{\tau}{\rightarrow} A \overset{g}{\rightarrow} B$ is the set
		\begin{align*}
			A^*_{B,R} :=A^*:= \left\{x \in B \mid \Delta(x) \in \overline{\ker \varphi}\right\}=\Delta^{-1}(\overline{\ker\varphi}).
		\end{align*}
	\end{defi}
	If $A^* = g(A)$ then $A$ is said to be saturated in $B$. Let us review some fundamental facts regarding relative Lipschitz saturation..
	
	\begin{prop}\cite{lipman1}\label{prop_satsubanel}
		$A^*_{B,R}$ is an $R$-subalgebra of $B$ which contains $g(A)$.
	\end{prop}

	\begin{prop}\cite{lipman1}\label{202306222352}
		Let $A$ be an $R$-subalgebra of $B$ and let $C$ be an $R$-subalgebra of $B$ containing $A$ as an $R$-subalgebra. \[\begin{tikzcd}
			R & A & C & B
			\arrow["\tau", from=1-1, to=1-2]
			\arrow[hook, from=1-2, to=1-3]
			\arrow[hook, from=1-3, to=1-4]
			\arrow["g"', curve={height=12pt}, hook, from=1-2, to=1-4]
		\end{tikzcd}\]

		If $C^*=C^*_{B,R}$ then $A^* \subseteq C^*$.
	\end{prop}

	\begin{cor}\cite{lipman1}\label{202306230013}
		If $A$ is an $R$-subalgebra of $B$ then: 
		
		\begin{enumerate}
			\item $A \subseteq A^*$.
			
			\item $(A^*)^* = A^*$.
		\end{enumerate}
	\end{cor}

	Let us establish some notation, which we will use throughout this work. Suppose that \[\begin{tikzcd}
		R & A & B \\
		{R'} & {A'} & {B'}
		\arrow["\tau", from=1-1, to=1-2]
		\arrow["g", from=1-2, to=1-3]
		\arrow["{\tau'}", from=2-1, to=2-2]
		\arrow["{g'}", from=2-2, to=2-3]
		\arrow["{f_R}"', from=1-1, to=2-1]
		\arrow["{f}", from=1-3, to=2-3]
		\arrow["{f_A}", from=1-2, to=2-2]
	\end{tikzcd}\]
	
	\noindent is a commutative diagram of ring morphisms. Then we can consider the diagram 
	 % https://q.uiver.app/#q=WzAsMTEsWzAsMCwiUiJdLFsxLDAsIkEiXSxbMiwwLCJCIl0sWzQsMCwiQlxcdW5kZXJzZXR7Un17XFxvdGltZXN9QiJdLFs2LDAsIkJcXHVuZGVyc2V0e0F9e1xcb3RpbWVzfUIiXSxbMCwyLCJSJyJdLFsxLDIsIkEnIl0sWzIsMiwiQiciXSxbNCwyLCJCJ1xcdW5kZXJzZXR7Uid9e1xcb3RpbWVzfUInIl0sWzYsMiwiQidcXHVuZGVyc2V0e0EnfXtcXG90aW1lc31CJyJdLFs0LDEsIkInXFx1bmRlcnNldHtSfXtcXG90aW1lc31CJyJdLFswLDEsIlxcdGF1Il0sWzEsMiwiZyJdLFsyLDMsIlxcRGVsdGEiXSxbMyw0LCJcXHZhcnBoaSJdLFswLDUsImZfUiIsMl0sWzEsNiwiZl9BIiwyXSxbMiw3LCJmIl0sWzUsNiwiXFx0YXUnIiwyXSxbNiw3LCJnJyIsMl0sWzcsOCwiXFxEZWx0YSciLDJdLFszLDEwLCJmXFx1bmRlcnNldHtSfXtcXG90aW1lc31mIl0sWzEwLDgsInAiXSxbMyw4LCJcXGJhcntmfSIsMix7ImN1cnZlIjo1fV0sWzgsOSwiXFx2YXJwaGknIl1d
	\[\begin{tikzcd}
		R & A & B && {B\underset{R}{\otimes}B} && {B\underset{A}{\otimes}B} \\
		&&&& {B'\underset{R}{\otimes}B'} \\
		{R'} & {A'} & {B'} && {B'\underset{R'}{\otimes}B'} && {B'\underset{A'}{\otimes}B'}
		\arrow["\tau", from=1-1, to=1-2]
		\arrow["g", from=1-2, to=1-3]
		\arrow["\Delta", from=1-3, to=1-5]
		\arrow["\varphi", from=1-5, to=1-7]
		\arrow["{f_R}"', from=1-1, to=3-1]
		\arrow["{f_A}"', from=1-2, to=3-2]
		\arrow["f", from=1-3, to=3-3]
		\arrow["{\tau'}"', from=3-1, to=3-2]
		\arrow["{g'}"', from=3-2, to=3-3]
		\arrow["{\Delta'}"', from=3-3, to=3-5]
		\arrow["{f\underset{R}{\otimes}f}", from=1-5, to=2-5]
		\arrow["p", from=2-5, to=3-5]
		\arrow["{\bar{f}}"', curve={height=30pt}, from=1-5, to=3-5]
		\arrow["{\varphi'}", from=3-5, to=3-7]
	\end{tikzcd}\eqno{(\clubsuit)}\]

\noindent  where $\varphi,\varphi'$ and $p$ are the canonical morphisms, and $\bar{f}:=p\circ(f\ten_Rf)$. Notice that for all $b\in B$ one has $$\bar{f}\circ\Delta(b)=\bar{f}(b\ten_R1_B-1_B\ten_Rb)=p(f(b)\ten_R1_{B'}-1_{B'}\ten_Rf(b))$$$$=f(b)\ten_{R'}1_{B'}-1_{B'}\ten_{R'}f(b)=\Delta'(f(b)).$$

Therefore, $(\clubsuit)$ is a commutative diagram. Lipman proved the following proposition.
	
	\begin{prop}\cite{lipman1} \label{prop3}
		Suppose that
		\[\begin{tikzcd}
			R & A & B \\
			{R'} & {A'} & {B'}
			\arrow["\tau", from=1-1, to=1-2]
			\arrow["g", from=1-2, to=1-3]
			\arrow["{\tau'}", from=2-1, to=2-2]
			\arrow["{g'}", from=2-2, to=2-3]
			\arrow["{f_R}"', from=1-1, to=2-1]
			\arrow["{f}", from=1-3, to=2-3]
			\arrow["{f_A}", from=1-2, to=2-2]
		\end{tikzcd}\]
		
		\noindent is a commutative diagram of ring morphisms. Then
		$$f\prtt{A^*_{B,R}} \subseteq \prtt{A'}^*_{B',R'}.$$
	\end{prop}

	Here we point out that Proposition \ref{202306222352} proved by Lipman is true even on the case where the sequence of rings $A\rightarrow C\rightarrow B$ is not necessarily a chain of subalgebras. Indeed, if % https://q.uiver.app/#q=WzAsNCxbMCwwLCJSIl0sWzEsMCwiQSJdLFsyLDAsIkMiXSxbMywwLCJCIl0sWzAsMSwiXFx0YXUiXSxbMSwyLCJcXGxhbWJkYSJdLFsyLDMsImdfQyJdLFsxLDMsImciLDIseyJjdXJ2ZSI6Mn1dXQ==
	$\begin{tikzcd}
		R & A & C & B
		\arrow["\tau", from=1-1, to=1-2]
		\arrow["\lambda", from=1-2, to=1-3]
		\arrow["{g_C}", from=1-3, to=1-4]
		\arrow["g"', curve={height=12pt}, from=1-2, to=1-4]
	\end{tikzcd}$ is a sequence of ring morphisms then we have the induced commutative diagram % https://q.uiver.app/#q=WzAsNixbMCwwLCJSIl0sWzEsMCwiQSJdLFsyLDAsIkMiXSxbMCwxLCJSIl0sWzEsMSwiQyJdLFsyLDEsIkIiXSxbMCwxLCJcXHRhdSJdLFsxLDIsIlxcbGFtYmRhIl0sWzIsNSwiZ19DIl0sWzAsMywiXFx0ZXh0cm17aWR9X1IiLDJdLFszLDQsIlxcdGF1X0MiLDJdLFsxLDQsIlxcbGFtYmRhIiwyXSxbNCw1LCJnX0MiLDJdXQ==
$\begin{tikzcd}
R & A & B \\
R & C & B
\arrow["\tau", from=1-1, to=1-2]
\arrow["g", from=1-2, to=1-3]
\arrow["{\textrm{id}_B}", from=1-3, to=2-3]
\arrow["{\textrm{id}_R}"', from=1-1, to=2-1]
\arrow["{\tau_C}"', from=2-1, to=2-2]
\arrow["\lambda" ', from=1-2, to=2-2]
\arrow["{g_C}"', from=2-2, to=2-3]
\end{tikzcd}$. Now it suffices to apply Proposition \ref{prop3} to get $A^*_{B,R}\sub C^*_{B,R}$.
	
	In order to investigate the contraction property, it is important to notice that Proposition \ref{prop3} says that $A^*_{B,R}\sub f^{-1}((A')^*_{B',R'})$, and its proof relies upon the fact that $\ker\varphi(B'\ten_{R'}B')\sub\ker\varphi'$,\footnote{For a ring morphism $h:S\rightarrow T$ and an ideal $I$ of $S$, we use the notation $IT$ to mean the ideal of $T$ generated by $h(I)$.} which is a consequence of how the kernels are generated. Indeed, we know 
	
	\begin{itemize}
		\item $\ker\varphi$ is generated by $\{\Delta(g(a))\mid a\in A\}$;
		
		\item $\ker\varphi'$ is generated by $\{\Delta'(g'(a'))\mid a'\in A'\}$,
	\end{itemize}

\noindent and since ($\clubsuit$) is a commutative diagram, for every $a\in A$ one has: 

$$\bar{f}(\Delta(g(a)))=\Delta'(f(g(a)))=\Delta'(g'(f_A(a)))\in \ker\varphi'.$$ 

\noindent This proves that $\bar{f}(\ker\varphi)\sub \ker\varphi'$, and consequently, $\ker\varphi(B'\ten_{R'}B')\sub\ker\varphi'$.
	
	The main goal of this work is to present sufficient conditions in order to get the \textit{contraction property} $$f^{-1}((A')^*_{B',R'})=A^*_{B,R},$$
	
	\noindent which will generalize the contractions obtained by Lipman in \cite{lipman1}. One of these conditions is $\ker\varphi(B'\ten_{R'}B')$ to be a reduction of $\ker\varphi'$, i.e, $\ker\varphi(B'\ten_{R'}B')\sub\ker\varphi'\sub\overline{\ker\varphi(B'\ten_{R'}B')}$, which is equivalent to $$\overline{\ker\varphi'}=\overline{\ker\varphi(B'\ten_{R'}B')}.$$
	
	\begin{defi}[Lipman diagram]\label{202306211340}
		We say that % https://q.uiver.app/#q=WzAsNixbMCwwLCJSIl0sWzEsMCwiQSJdLFsyLDAsIkIiXSxbMCwxLCJSJyJdLFsxLDEsIkEnIl0sWzIsMSwiQiciXSxbMCwxLCJcXHRhdSJdLFsxLDIsImciXSxbMyw0LCJcXHRhdSciLDJdLFs0LDUsImcnIiwyXSxbMCwzLCJmX1IiLDJdLFsxLDQsImZfQSIsMl0sWzIsNSwiZiIsMl1d
		\[\begin{tikzcd}
			R & A & B \\
			{R'} & {A'} & {B'}
			\arrow["\tau", from=1-1, to=1-2]
			\arrow["g", from=1-2, to=1-3]
			\arrow["{\tau'}", from=2-1, to=2-2]
			\arrow["{g'}", from=2-2, to=2-3]
			\arrow["{f_R}"', from=1-1, to=2-1]
			\arrow["{f}", from=1-3, to=2-3]
			\arrow["{f_A}", from=1-2, to=2-2]
		\end{tikzcd}\] is a \textbf{Lipman diagram} if it is a commutative diagram of ring morphisms and $$\overline{\ker\varphi'}=\overline{\ker\varphi(B'\ten_{R'}B')}.$$
		
		\noindent	We say that it is a \textbf{strong Lipman diagram} if $\ker\varphi'=\ker\varphi(B'\ten_{R'}B').$
	\end{defi} 
	
	In  next result, we give a sufficient condition on $f_A$ to get strong Lipman diagrams, which is useful to generate some examples for cases where $A'\neq A$.
	
	\begin{prop}\label{202306011826}
		Consider the commutative diagram % https://q.uiver.app/#q=WzAsNixbMCwwLCJSIl0sWzEsMCwiQSJdLFsyLDAsIkIiXSxbMCwxLCJSJyJdLFsxLDEsIkEnIl0sWzIsMSwiQiciXSxbMCwxLCJcXHRhdSJdLFsxLDIsImciXSxbMyw0LCJcXHRhdSciLDJdLFs0LDUsImcnIiwyXSxbMCwzLCJmX1IiLDJdLFsxLDQsImZfQSIsMl0sWzIsNSwiZiIsMl1d
		\[\begin{tikzcd}
			R & A & B \\
			{R'} & {A'} & {B'}
			\arrow["\tau", from=1-1, to=1-2]
			\arrow["g", from=1-2, to=1-3]
			\arrow["{\tau'}", from=2-1, to=2-2]
			\arrow["{g'}", from=2-2, to=2-3]
			\arrow["{f_R}"', from=1-1, to=2-1]
			\arrow["{f}", from=1-3, to=2-3]
			\arrow["{f_A}", from=1-2, to=2-2]
		\end{tikzcd} \eqno{(\star)}\] \noindent of ring morphisms. If $f_A$ is surjective then $(\star)$ is a strong Lipman diagram.
	\end{prop}
	
	\begin{proof}
		We only have to check that $\ker\varphi'\sub\ker\varphi(B'\ten_RB')$.	We already know that: 
		
		\begin{itemize}
			\item $\ker\varphi$ is generated by $\{\Delta(g(a))\mid a\in A\}$;
			
			\item $\ker\varphi'$ is generated by $\{\Delta'(g'(a'))\mid a'\in A'\}$. 
		\end{itemize}
		
		Let $a'\in A'$. Since $f_A$ is surjective then $a'=f_A(a)$ for some $a\in A$. Thus: $$\Delta'(g'(a))=\Delta'(g'(f(a)))=\Delta'(f(g(a)))=(\Delta'\circ f)(g(a))$$$$=\bar{f}(\Delta(g(a)))\in\bar{f}(\ker\varphi)\sub \ker\varphi(B'\ten_{R'}B').$$	
	\end{proof}
	
	In \cite{lipman1}, Lipman studied a condition for which a sequence of ring morphism % https://q.uiver.app/#q=WzAsNCxbMCwwLCJSIl0sWzEsMCwiQSJdLFsyLDAsIkIiXSxbMywwLCJCJyJdLFswLDEsIlxcdGF1Il0sWzEsMiwiZyJdLFsyLDMsImYiXV0=
	\[\begin{tikzcd}
		R & A & B & {B'}
		\arrow["\tau", from=1-1, to=1-2]
		\arrow["g", from=1-2, to=1-3]
		\arrow["f", from=1-3, to=1-4]
	\end{tikzcd}\] \noindent satisfies $f^{-1}(A^*_{B',R})=A^*_{B,R}$. The first step to obtaining our generalization for some of the Lipman results is to show that every such sequence can be seen as a Lipman diagram. 
	
	\begin{cor}\label{202306011657}
		If % https://q.uiver.app/#q=WzAsNCxbMCwwLCJSIl0sWzEsMCwiQSJdLFsyLDAsIkIiXSxbMywwLCJCJyJdLFswLDEsIlxcdGF1Il0sWzEsMiwiZyJdLFsyLDMsImYiXV0=
		\begin{tikzcd}
			R & A & B & {B'}
			\arrow["\tau", from=1-1, to=1-2]
			\arrow["g", from=1-2, to=1-3]
			\arrow["f", from=1-3, to=1-4]
		\end{tikzcd} is a sequence of ring morphism then % https://q.uiver.app/#q=WzAsNixbMCwwLCJSIl0sWzEsMCwiQSJdLFsyLDAsIkIiXSxbMCwxLCJSIl0sWzEsMSwiQSJdLFsyLDEsIkInIl0sWzAsMSwiXFx0YXUiXSxbMSwyLCJnIl0sWzMsNCwiXFx0YXUiLDJdLFs0LDUsImZcXGNpcmMgZyIsMl0sWzAsMywiXFx0ZXh0cm17aWR9X1IiLDJdLFsxLDQsIlxcdGV4dHJte2lkfV9BIiwyXSxbMiw1LCJmIiwyXV0=
		\[\begin{tikzcd}
			R & A & B \\
			R & A & {B'}
			\arrow["\tau", from=1-1, to=1-2]
			\arrow["g", from=1-2, to=1-3]
			\arrow["\tau"', from=2-1, to=2-2]
			\arrow["{f\circ g}"', from=2-2, to=2-3]
			\arrow["{\textrm{id}_R}"', from=1-1, to=2-1]
			\arrow["{\textrm{id}_A}"', from=1-2, to=2-2]
			\arrow["f"', from=1-3, to=2-3]
		\end{tikzcd}\] is a strong Lipman diagram.
	\end{cor}
	
	\begin{proof} 
		It is an immediate consequence of Proposition \ref{202306011826}, once $\id_A$ is surjective.
	\end{proof}
	
	The next example was introduced by Lipman in \cite{lipman1} in a context where he was interested in understanding when a change of base would preserve Lipschitz saturation. In our terminology, this example is a (strong) Lipman diagram.
	
	\begin{example}\normalfont
		Let $R\overset{\tau}{\rightarrow}A\overset{g}{\rightarrow}B$ be a sequence of ring morphisms, and let $R'$ be an $R$-algebra with structure morphism $\delta_R:R\rightarrow R'$. Denote $A':=A\ten_RR'$, $B':=B\ten_RR'$ and consider the canonical morphisms $\delta_A:A\rightarrow A'$ and $\delta:B\rightarrow B'$. Besides, consider $g':A'\rightarrow B'$ given by $g':=g\ten_R\id_{R'}$, and let $\tau':R'\rightarrow A'$ be the morphism given by $\tau'(r'):=1_A\ten_Rr', \forall r'\in R'$. Notice that for every $r\in R$ we have $$\delta_A\circ\tau(r)=\tau(r)\ten_R1_{R'}=(r1_A)\ten_R1_{R'}=1_A\ten_R(r1_{R'})=1_A\ten_R\delta_R(r)=\tau'\circ\delta_R(r).$$
		
		Thus, $\delta_A\circ\tau=\tau'\circ\delta_R$. Further, for every $a\in A$ one has $$\delta\circ g(a)=g(a)\ten_R1_{R'}=(g\ten_R\id_{R'})(a\ten_R1_{R'})=g'\circ\delta_A(a).$$
		
		Hence, \[\begin{tikzcd}
			R & A & B \\
			{R'} & {A'} & {B'}
			\arrow["\tau", from=1-1, to=1-2]
			\arrow["g", from=1-2, to=1-3]
			\arrow["{\tau'}", from=2-1, to=2-2]
			\arrow["{g'}", from=2-2, to=2-3]
			\arrow["{\delta_R}"', from=1-1, to=2-1]
			\arrow["{\delta}", from=1-3, to=2-3]
			\arrow["{\delta_A}", from=1-2, to=2-2]
		\end{tikzcd} \eqno{(\star)}\]
		
		\noindent is a commutative diagram, and we can consider the induced diagram $(\clubsuit)$ for this case. 
		
		Let us prove that $(\star)$ is a strong Lipman diagram. Indeed we know that $\ker\varphi$ is generated by $\Delta(g(A))$ and $\ker\varphi'$ is generated by $\Delta'(g'(A'))$. Since $A'=A\ten_RR'$ then $\ker\varphi'$ is generated by $\{a\ten_Rr'\mid a\in A\mbox{ e }r'\in R'\}$. Thus, for all $a\in A$ and $r'\in R'$ one has $$\Delta'(g'(a\ten_Rr'))=\Delta'(g(a)\ten_Rr')=\Delta'((1_A\ten_Rr')(g(a)\ten_R1_{R'}))=\Delta'((\tau'(r'))\delta(g(a)))$$$$=\Delta'(r'\cdot \delta(g(a)))=r'\cdot \Delta'(\delta(g(a)))=r'\cdot \bar{\delta}(\Delta(g(a)))\in\bar{\delta}(\ker\varphi).$$
		
		Therefore, $\ker\varphi'=\ker\varphi(B'\ten_{R'}B')$ and $(\star)$ is a strong Lipman diagram.		
	\end{example}

	\section{Universally injective contraction}
	
	In \cite{lipman1}, Lipman proved that if we have a sequence % https://q.uiver.app/#q=WzAsNCxbMCwwLCJSIl0sWzEsMCwiQSJdLFsyLDAsIkIiXSxbMywwLCJCJyJdLFswLDEsIlxcdGF1Il0sWzEsMiwiZyJdLFsyLDMsImYiXV0=
	\begin{tikzcd}
		R & A & B & {B'}
		\arrow["\tau", from=1-1, to=1-2]
		\arrow["g", from=1-2, to=1-3]
		\arrow["f", from=1-3, to=1-4]
	\end{tikzcd} of ring morphisms and if $f$ is faithfully flat then $f^{-1}((A')^*_{B',R'})=A^*_{B,R}$. 
	
	In this section, we use the universally injective condition (more general than faithfully flatness) as an important condition so that Lipman diagrams contracts. 
	
	Recall that a ring morphism $\alpha:S\rightarrow S'$ is universally injective if and only if for every $S$-module $N$ the induced morphism

	$$\begin{matrix}
		\tilde{\alpha}: & N & \longrightarrow  &  N\ten_SS'\\
		& n & \longmapsto         & n\ten_S 1_{S'}
	\end{matrix}$$

	\noindent is injective.

	In this way, we need some basic results on commutative algebra. For some of them, we were not able to find a satisfactory reference, and because of this, we present our own proof for them.
	
	\begin{lema}
		Let $\alpha:S\rightarrow S'$ be a morphism of $R$-algebras and let $T$ be an $R$-algebra. If $\alpha$ is universally injective then $\alpha\ten_R\id_T:S\ten_RT\rightarrow S'\ten_RT$ is universally injective.
	\end{lema}
	
	\begin{proof}
		Let $N$ be a $S\ten_RT$-module.  The canonical morphism $S\rightarrow S\ten_RT$ induces a $S$-module structure on $N$, and by hypothesis $\tilde{\alpha}:N\rightarrow N\ten_SS'$ is injective. Notice that $$N\underset{S\ten_RT}{\otimes}(S'\ten_RT)\cong N\underset{S\ten_RT}{\otimes}(S'\ten_S(S\ten_RT))\cong(N\underset{S\ten_RT}{\otimes}(S\ten_RT))\ten_SS'\cong N\ten_SS',$$
		
		\noindent and this isomorphism is such that 
		
		% https://q.uiver.app/#q=WzAsMyxbMCwxLCJOIl0sWzIsMCwiTlxcdW5kZXJzZXR7U1xcdW5kZXJzZXR7Un17XFxvdGltZXN9VH17XFxvdGltZXN9KFMnXFx1bmRlcnNldHtSfXtcXG90aW1lc31UKSJdLFsyLDEsIk5cXHVuZGVyc2V0e1N9e1xcb3RpbWVzfVMnIl0sWzAsMSwiXFx3aWRldGlsZGV7XFxhbHBoYVxcdW5kZXJzZXR7Un17XFxvdGltZXN9XFx0ZXh0cm17aWR9X1R9IiwwLHsiY3VydmUiOi0yfV0sWzAsMiwiXFx0aWxkZXtcXGFscGhhfSIsMl0sWzEsMiwiXFxjb25nIl1d
		\[\begin{tikzcd}
			&& {N\underset{S\underset{R}{\otimes}T}{\otimes}(S'\underset{R}{\otimes}T)} \\
			N && {N\underset{S}{\otimes}S'}
			\arrow["{\widetilde{\alpha\underset{R}{\otimes}\textrm{id}_T}}", curve={height=-12pt}, from=2-1, to=1-3]
			\arrow["{\tilde{\alpha}}"', from=2-1, to=2-3]
			\arrow["\cong", from=1-3, to=2-3]
		\end{tikzcd}\]

		\noindent is commutative. Hence, $\widetilde{\alpha\underset{R}{\otimes}\textrm{id}_T}$ is injective, as desired.
	\end{proof}
	
	\begin{cor}\label{202305231737}
		If $\alpha:S\rightarrow S'$ and $\beta:T\rightarrow T'$ are universally injective morphisms of $R$-algebras then $$\alpha\ten_R\beta:S\ten_RT\rightarrow S'\ten_RT'$$
		
		\noindent is universally injective.
	\end{cor}
	
	\begin{proof}
		Since universal injectiveness is preserved by the composition of maps, this corollary is a consequence of the previous lemma, once % https://q.uiver.app/#q=WzAsMyxbMCwwLCJTXFx1bmRlcnNldHtSfXtcXG90aW1lc31UIl0sWzQsMCwiUydcXHVuZGVyc2V0e1J9e1xcb3RpbWVzfVQnIl0sWzIsMiwiUydcXHVuZGVyc2V0e1J9e1xcb3RpbWVzfVQiXSxbMCwxLCJcXGFscGhhXFx1bmRlcnNldHtSfXtcXG90aW1lc31cXGJldGEiXSxbMCwyLCJcXGFscGhhXFx1bmRlcnNldHtSfXtcXG90aW1lc31cXHRleHRybXtpZH1fVCIsMl0sWzIsMSwiXFx0ZXh0cm17aWR9X3tTJ31cXHVuZGVyc2V0e1J9e1xcb3RpbWVzfVxcYmV0YSIsMl1d
		\[\begin{tikzcd}
			{S\underset{R}{\otimes}T} &&&& {S'\underset{R}{\otimes}T'} \\
			\\
			&& {S'\underset{R}{\otimes}T}
			\arrow["{\alpha\underset{R}{\otimes}\beta}", from=1-1, to=1-5]
			\arrow["{\alpha\underset{R}{\otimes}\textrm{id}_T}"', from=1-1, to=3-3]
			\arrow["{\textrm{id}_{S'}\underset{R}{\otimes}\beta}"', from=3-3, to=1-5]
		\end{tikzcd}\] \noindent is commutative. 
	\end{proof}
	
	%\begin{lema}
	%If $\alpha:S\rightarrow S'$ is a flat (faithfully flat) ring morphism and $T$ is $S'$-algebra then the canonical morphism $p:T\ten_ST\rightarrow T\ten_{S'}T$ is flat (faithfully flat).
	%\end{lema}
	
	%\begin{proof}
	%	\DD
	%\end{proof}
	
	\begin{lema}\label{202306082227}
		Let $h: S\rightarrow T$ be a universally injective ring morphism. If $I$ is an ideal of $S$ then $$h^{-1}(\overline{IT})=\overline{I}.$$
	\end{lema}
	
	\begin{proof}
		First we show that $h^{-1}(JT)=J$, for all ideal $J$ of $S$. It is clear that $J\sub h^{-1}(JT)$. Conversely, consider the morphism $\bar{h}:\dfrac{S}{J}\rightarrow \dfrac{T}{JT}$ induced by $h$. Since $h$ is universally injective and $\dfrac{S}{J}$ is a $S$-module then $$\begin{matrix}
			\tilde{h}: & \dfrac{S}{J} & \longrightarrow & \dfrac{S}{J}\ten_ST\\
			& u & \longmapsto & u\ten_S1_T
		\end{matrix}$$\noindent is injective. The canonical isomorphism $\phi:\dfrac{S}{J}\ten_ST\rightarrow \dfrac{T}{JT}$ clearly makes the diagram % https://q.uiver.app/#q=WzAsMyxbMCwwLCJcXGRmcmFje1N9e0l9Il0sWzIsMCwiXFxkZnJhY3tUfXtJVH0iXSxbMSwxLCJcXGRmcmFje1N9e0l9XFxvdGltZXNfU1QiXSxbMCwxLCJcXGJhcntofSJdLFswLDIsIlxcdGlsZGV7aH0iLDJdLFsyLDEsIlxccGhpIiwyXV0=
		\[\begin{tikzcd}
			{\dfrac{S}{J}} && {\dfrac{T}{JT}} \\
			& {\dfrac{S}{J}\otimes_ST}
			\arrow["{\bar{h}}", from=1-1, to=1-3]
			\arrow["{\tilde{h}}"', from=1-1, to=2-2]
			\arrow["\phi"', from=2-2, to=1-3]
		\end{tikzcd}\] \noindent to be commutative. Thus, $\bar{h}$ is injective, which implies that $h^{-1}(JT)\sub J$. 
		
		Now, let us prove that $h^{-1}(\overline{IT})=\overline{I}$. The persistence of the integral closure (see 1.1.3 in \cite{sh}) gives the inclusion $\overline{I}\sub h^{-1}(\overline{IT})$. Conversely, let $x\in h^{-1}(\overline{IT})$. Thus $h(x)\in \overline{IT}$, and an integral dependence equation of $h(x)$ over $IT$ implies there is $n\in\bN$ such that $$(h(x))^n\in I(I+(x))^{n-1}T,$$ \noindent i.e, $x^n\in h^{-1}(I(I+(x))^{n-1}T)$. Since $h$ is universally injective, we just saw that\\ $h^{-1}(I(I+(x))^{n-1}T)=I(I+(x))^{n-1}$, so $x^n\in I(I+(x))^{n-1}$. Therefore (see Prop. 1.1.7 in \cite{sh}), $x\in \overline{I}$.   
	\end{proof}
	
	\begin{lema}\label{202306082225}
		If $f_R: R\rightarrow R'$ is a universally injective ring morphism and $B'$ is an $R'$-algebra then the canonical morphism $p: B'\ten_RB'\rightarrow B'\ten_{R'}B'$ is an isomorphism.
	\end{lema}
	
	\begin{proof}
		Once $p$ is surjective, we only have to show that $p$ is injective. Since $B'\ten_RB'$ is an $R$-module and $f_R$ is universally injective then $\tilde{f_R}:B'\ten_RB'\rightarrow (B'\ten_RB')\ten_RR'$ is injective. Besides, by \ref{202305231737} we know that $f_R\ten_Rf_R:R\ten_RR\cong R\rightarrow R'\ten_RR'$ is universally injective, and since $B'\ten_{R'}B'$ is an $R$-module then $$\widetilde{f_R\ten_Rf_R}: B'\ten_{R'}B'\rightarrow (B'\ten_{R'}B')\ten_R(R'\ten_RR')$$
		
		\noindent is injective. There exists a canonical isomorphism $\phi:(B'\ten_RB')\ten_RR'\rightarrow (B'\ten_{R'}B')\ten_R(R'\ten_RR')$ such that the diagram % https://q.uiver.app/#q=WzAsNCxbMCwwLCJCJ1xcdW5kZXJzZXR7Un17XFxvdGltZXN9QiciXSxbMiwwLCIoQidcXHVuZGVyc2V0e1J9e1xcb3RpbWVzfUInKVxcdW5kZXJzZXR7Un17XFxvdGltZXN9UiciXSxbMCwxLCJCJ1xcdW5kZXJzZXR7Uid9e1xcb3RpbWVzfUInIl0sWzIsMSwiKEInXFx1bmRlcnNldHtSJ317XFxvdGltZXN9QicpXFx1bmRlcnNldHtSfXtcXG90aW1lc30oUidcXHVuZGVyc2V0e1J9e1xcb3RpbWVzfVInKSJdLFswLDEsIlxcdGlsZGV7Zl9SfSJdLFswLDIsInAiLDJdLFsyLDMsIlxcd2lkZXRpbGRle2ZfUlxcdW5kZXJzZXR7Un17XFxvdGltZXN9Zl9SfSJdLFsxLDMsIlxccGhpIl1d
		\[\begin{tikzcd}
			{B'\underset{R}{\otimes}B'} && {(B'\underset{R}{\otimes}B')\underset{R}{\otimes}R'} \\
			{B'\underset{R'}{\otimes}B'} && {(B'\underset{R'}{\otimes}B')\underset{R}{\otimes}(R'\underset{R}{\otimes}R')}
			\arrow["{\tilde{f_R}}", from=1-1, to=1-3]
			\arrow["p"', from=1-1, to=2-1]
			\arrow["{\widetilde{f_R\underset{R}{\otimes}f_R}}", from=2-1, to=2-3]
			\arrow["\phi", from=1-3, to=2-3]
		\end{tikzcd}\]\noindent commutes. Therefore, $p$ is injective.
	\end{proof}

	Now we state the main result of this section.
	
	\begin{teo}[Universally injective contraction]\label{202306011654} 	Suppose that
		\[\begin{tikzcd}
			R & A & B \\
			{R'} & {A'} & {B'}
			\arrow["\tau", from=1-1, to=1-2]
			\arrow["g", from=1-2, to=1-3]
			\arrow["{\tau'}", from=2-1, to=2-2]
			\arrow["{g'}", from=2-2, to=2-3]
			\arrow["{f_R}"', from=1-1, to=2-1]
			\arrow["{f}", from=1-3, to=2-3]
			\arrow["{f_A}", from=1-2, to=2-2]
		\end{tikzcd}\]
		
		\noindent is a Lipman diagram. If $f_R$ and $f$ are universally injective then $$f^{-1}((A')^*_{B',R'})=A^*_{B,R}.$$
	\end{teo}
	
	\begin{proof}
		Let $\bar{f}$ be the composition of the morphisms % https://q.uiver.app/#q=WzAsMyxbMCwwLCJCXFx1bmRlcnNldHtSfXtcXG90aW1lc31CIl0sWzIsMCwiQidcXHVuZGVyc2V0e1J9e1xcb3RpbWVzfUInIl0sWzQsMCwiQidcXHVuZGVyc2V0e1InfXtcXG90aW1lc31CJyJdLFswLDEsImZcXHVuZGVyc2V0e1J9e1xcb3RpbWVzfWYiXSxbMSwyLCJwIl1d
		\begin{tikzcd}
			{B\underset{R}{\otimes}B} & {B'\underset{R}{\otimes}B'} & {B'\underset{R'}{\otimes}B'}
			\arrow["{f\underset{R}{\otimes}f}", from=1-1, to=1-2]
			\arrow["p", from=1-2, to=1-3]
		\end{tikzcd}. By Corollary \ref{202305231737} we have that $f\ten_R f$ is universally injective. Besides, Lemma \ref{202306082225} implies that $p$ is an isomorphism. So, $\bar{f}=p\circ(f\ten_Rf)$ is also universally injective. Thus, Lemma \ref{202306082227} ensures that $$\bar{f}^{-1}(\overline{\ker\varphi'})=\bar{f}^{-1}(\overline{\ker\varphi(B'\ten_{R}B')})=\overline{\ker\varphi}.$$
		
		It is easy to see that % https://q.uiver.app/#q=WzAsNCxbMiwwLCJCXFx1bmRlcnNldHtSfXtcXG90aW1lc31CIl0sWzIsMSwiQidcXHVuZGVyc2V0e1InfXtcXG90aW1lc31CJyJdLFswLDAsIkIiXSxbMCwxLCJCJyJdLFsyLDAsIlxcRGVsdGEiXSxbMywxLCJcXERlbHRhJyJdLFsyLDMsImYiLDJdLFswLDEsIlxcYmFye2Z9Il1d
		\[\begin{tikzcd}
			B && {B\underset{R}{\otimes}B} \\
			{B'} && {B'\underset{R'}{\otimes}B'}
			\arrow["\Delta", from=1-1, to=1-3]
			\arrow["{\Delta'}", from=2-1, to=2-3]
			\arrow["f"', from=1-1, to=2-1]
			\arrow["{\bar{f}}", from=1-3, to=2-3]
		\end{tikzcd}\] \noindent commutes, and consequently
		
		$$f^{-1}((A')^*_{B', R'})=f^{-1}(\Delta'^{-1}(\overline{\ker\varphi'}))=\Delta^{-1}(\bar{f}^{-1}(\overline{\ker\varphi'}))=\Delta^{-1}(\overline{\ker\varphi})=A^*_{B,R}.$$
	\end{proof}
	
	As a consequence, we obtain a generalization of the contraction for the particular case studied by Lipman, since it is very known that every faithfully flat morphism is universally injective, although it is not equivalent conditions.
	
	\begin{cor}\label{202305241951}
		Let $R\overset{\tau}{\rightarrow} A \overset{g}{\rightarrow}B\overset{f}{\rightarrow}B'$ be a sequence of ring morphisms, and assume that $f$ is universally injective. Then: $$A^*_{B,R}=f^{-1}(A^*_{B',R}).$$
	\end{cor}
	
	\begin{proof}
		In Corollary \ref{202306011657} we have seen that the above sequence gives rise to a Lipman diagram % https://q.uiver.app/#q=WzAsNixbMCwwLCJSIl0sWzEsMCwiQSJdLFsyLDAsIkIiXSxbMCwxLCJSIl0sWzEsMSwiQSJdLFsyLDEsIkInIl0sWzAsMSwiXFx0YXUiXSxbMSwyLCJnIl0sWzMsNCwiXFx0YXUiLDJdLFs0LDUsImZcXGNpcmMgZyIsMl0sWzAsMywiXFx0ZXh0cm17aWR9X1IiLDJdLFsxLDQsIlxcdGV4dHJte2lkfV9BIiwyXSxbMiw1LCJmIiwyXV0=
		\[\begin{tikzcd}
			R & A & B \\
			R & A & {B'}
			\arrow["\tau", from=1-1, to=1-2]
			\arrow["g", from=1-2, to=1-3]
			\arrow["\tau"', from=2-1, to=2-2]
			\arrow["{f\circ g}"', from=2-2, to=2-3]
			\arrow["{\textrm{id}_R}"', from=1-1, to=2-1]
			\arrow["{\textrm{id}_A}"', from=1-2, to=2-2]
			\arrow["f"', from=1-3, to=2-3]
		\end{tikzcd}\].
		
		In this case, $f_R=\id_R$ and $f$ are universally injective, and the proof is done.
	\end{proof}

	Now we apply the main theorem of this section to get a new class of algebras where the conclusion of Proposition 1.3 of \cite{lipman1} is satisfied.
	
	\begin{prop}
		Consider a sequence of ring morphisms % https://q.uiver.app/#q=WzAsNCxbMCwwLCJSIl0sWzEsMCwiUiciXSxbMiwwLCJBIl0sWzMsMCwiQiJdLFswLDEsImZfUiJdLFsxLDJdLFsyLDMsImciXV0=
		$\begin{tikzcd}
			R & {R'} & A & B
			\arrow["{f_R}", from=1-1, to=1-2]
			\arrow[from=1-2, to=1-3]
			\arrow["g", from=1-3, to=1-4]
		\end{tikzcd}$. If $f_R$ is universally injective then $$A^*_{B,R}=A^*_{B,R'}.$$ 
	\end{prop}

\begin{proof}
	From the above sequence we get a strong Lipman diagram % https://q.uiver.app/#q=WzAsNixbMCwwLCJSIl0sWzIsMCwiQiJdLFswLDEsIlInIl0sWzEsMSwiQSJdLFsyLDEsIkIiXSxbMSwwLCJBIl0sWzAsNV0sWzAsMiwiZl9SIiwyXSxbMiwzXSxbMyw0LCJnIiwyXSxbNSwxLCJnIl0sWzUsMywiXFx0ZXh0cm17aWR9X0EiXSxbMSw0LCJcXHRleHRybXtpZH1fQiJdXQ==
	\[\begin{tikzcd}
		R & A & B \\
		{R'} & A & B
		\arrow[from=1-1, to=1-2]
		\arrow["{f_R}"', from=1-1, to=2-1]
		\arrow[from=2-1, to=2-2]
		\arrow["g"', from=2-2, to=2-3]
		\arrow["g", from=1-2, to=1-3]
		\arrow["{\textrm{id}_A}", from=1-2, to=2-2]
		\arrow["{\textrm{id}_B}", from=1-3, to=2-3]
	\end{tikzcd}\]\noindent  and since $\id_B$ is universally equivalent then Theorem \ref{202306011654} $\id_B^{-1}(A^*_{B,R'})=A^*_{B,R}$. Hence, $A^*_{B,R}=A^*_{B,R'}$.
\end{proof}

%In particular, if $R'$ is a faithfully flat $R$-algebra then $A^*_{B,R}=A^*_{B,R'}$.
	
	The next result follows directly from the Lemma \ref{202306082225} and Theorem \ref{202306011654} .

\begin{prop}
	Suppose that % https://q.uiver.app/#q=WzAsNSxbMCwwLCJSIl0sWzEsMCwiQSJdLFsyLDEsIkIiXSxbMCwyLCJSJyJdLFsxLDIsIkEnIl0sWzAsMywiZl9SIl0sWzEsNCwiZl9BIl0sWzEsMiwiZyJdLFs0LDIsImcnIiwyXSxbMyw0LCJcXHRhdSciXSxbMCwxLCJcXHRhdSJdXQ==
	\[\begin{tikzcd}
		R & A \\
		&& B \\
		{R'} & {A'}
		\arrow["{f_R}", from=1-1, to=3-1]
		\arrow["{f_A}", from=1-2, to=3-2]
		\arrow["g", from=1-2, to=2-3]
		\arrow["{g'}"', from=3-2, to=2-3]
		\arrow["{\tau'}", from=3-1, to=3-2]
		\arrow["\tau", from=1-1, to=1-2]
	\end{tikzcd} \eqno{(\star)}\] \noindent is a commutative diagram of ring morphisms. If $f_R$ and $f_A$ are universally injective then $$(A')^*_{B,R'}=A^*_{B,R}.$$
\end{prop}

\begin{proof}
	By Lemma \ref{202306082225} the canonical morphisms $p_R:B\ten_RB\rightarrow B\ten_{R'}B$ and $p_A:B\ten_AB\rightarrow B\ten_{A'}B$ are isomorphisms. Thus, the commutative diagram % https://q.uiver.app/#q=WzAsNCxbMCwwLCJCXFx1bmRlcnNldHtSfXtcXG90aW1lc31CIl0sWzEsMCwiQlxcdW5kZXJzZXR7QX17XFxvdGltZXN9QiJdLFswLDEsIkJcXHVuZGVyc2V0e1InfXtcXG90aW1lc31CIl0sWzEsMSwiQlxcdW5kZXJzZXR7QSd9e1xcb3RpbWVzfUIiXSxbMCwyLCJwX1IiLDJdLFsyLDMsIlxcdmFycGhpJyIsMl0sWzEsMywicF9BIl0sWzAsMSwiXFx2YXJwaGkiXV0=
	\[\begin{tikzcd}
		{B\underset{R}{\otimes}B} & {B\underset{A}{\otimes}B} \\
		{B\underset{R'}{\otimes}B} & {B\underset{A'}{\otimes}B}
		\arrow["{p_R}"', from=1-1, to=2-1]
		\arrow["{\varphi'}"', from=2-1, to=2-2]
		\arrow["{p_A}", from=1-2, to=2-2]
		\arrow["\varphi", from=1-1, to=1-2]
	\end{tikzcd}\]\noindent satisfies $p_R(\ker\varphi)=\ker\varphi'$. Hence, $(\star)$ is a strong Lipman diagram, and since $f_R$ and $\id_B$ are universally injective then Theorem \ref{202306011654} implies that $\id_B^{-1}((A')^*_{B,R'})=A^*_{B,R}$, i.e, $$(A')^*_{B,R'}=A^*_{B,R}.$$
\end{proof}

As a direct consequence, we obtain a version of the previous proposition formulated for the sequences analyzed by Lipman in \cite{lipman1}.

\begin{cor}
	If % https://q.uiver.app/#q=WzAsNCxbMCwwLCJSIl0sWzEsMCwiQSJdLFsyLDAsIkEnIl0sWzMsMCwiQiJdLFswLDFdLFsxLDIsImZfQSJdLFsyLDNdXQ==
	$\begin{tikzcd}
		R & A & {A'} & B
		\arrow[from=1-1, to=1-2]
		\arrow["{f_A}", from=1-2, to=1-3]
		\arrow[from=1-3, to=1-4]
	\end{tikzcd}$ \noindent is a sequence of ring morphisms and $f_A$ is universally injective then $$(A')^*_{B,R}=A^*_{B,R}.$$
\end{cor}

In Corollary \ref{202306230013} Lipman showed the idempotency $(A^*_{B,R})^*_{B,R}=A^*_{B,R}$ in the case where $A$ is an $R$-subalgebra of $B$. In the following, we employ Theorem \ref{202306011654} to extend the concept of idempotency to any sequence of rings.

\begin{prop}
	If % https://q.uiver.app/#q=WzAsMyxbMCwwLCJSIl0sWzEsMCwiQSJdLFsyLDAsIkIiXSxbMCwxLCJcXHRhdSJdLFsxLDIsImciXV0=
	$\begin{tikzcd}
		R & A & B
		\arrow["\tau", from=1-1, to=1-2]
		\arrow["g", from=1-2, to=1-3]
	\end{tikzcd}$ is a sequence of ring morphisms then $$(A^*_{B,R})^*_{B,R}=A^*_{B,R}.$$ 
\end{prop}

\begin{proof}
	Consider the ring morphism $f_A:A\rightarrow A^*_{B,R}$ given by $f_A(a):=g(a)$, $\forall  a\in A$ and let $g':A^*_{B,R}\hookrightarrow B$ be the inclusion map. Thus it is clear that % https://q.uiver.app/#q=WzAsNixbMCwwLCJSIl0sWzAsMSwiUiJdLFsxLDAsIkEiXSxbMSwxLCJBXipfe0IsUn0iXSxbMiwwLCJCIl0sWzIsMSwiQiJdLFswLDIsIlxcdGF1Il0sWzIsNCwiZyJdLFsxLDMsImZfQVxcY2lyY1xcdGF1IiwyXSxbMiwzLCJmX0EiLDJdLFswLDEsIlxcdGV4dHJte2lkfV9SIiwyXSxbNCw1LCJcXHRleHRybXtpZH1fQiJdLFszLDUsImcnIiwyLHsic3R5bGUiOnsidGFpbCI6eyJuYW1lIjoiaG9vayIsInNpZGUiOiJ0b3AifX19XV0=
	\[\begin{tikzcd}
		R & A & B \\
		R & {A^*_{B,R}} & B
		\arrow["\tau", from=1-1, to=1-2]
		\arrow["g", from=1-2, to=1-3]
		\arrow["{f_A\circ\tau}"', from=2-1, to=2-2]
		\arrow["{f_A}"', from=1-2, to=2-2]
		\arrow["{\textrm{id}_R}"', from=1-1, to=2-1]
		\arrow["{\textrm{id}_B}", from=1-3, to=2-3]
		\arrow["{g'}"', hook, from=2-2, to=2-3]
	\end{tikzcd}\eqno{(\star)}\]
	
	\noindent commutes. We know that the kernel of the morphism $\varphi':B\ten_RB\rightarrow B\ten_{A^*_{B,R}}B$ is generated by $\Delta(g'(A^*_{B,R}))=\Delta(A^*_{B,R})$. Since $A^*_{B,R}=\Delta^{-1}(\overline{\ker\varphi})$ then $\ker\varphi'\sub\overline{\ker\varphi}$. As a result, $(\star)$ becomes a Lipman diagram. Moreover, since $\id_R$ and $\id_B$ are universally injective then Theorem \ref{202306011654} ensures that $\id_B^{-1}((A^*_{B,R})^*_{B,R})=A^*_{B,R}$, implying $(A^*_{B,R})^*_{B,R}=A^*_{B,R}$.
\end{proof}

	\section{Integral contraction}
	
	In \cite{lipman1}, Lipman proved that if we have a sequence % https://q.uiver.app/#q=WzAsNCxbMCwwLCJSIl0sWzEsMCwiQSJdLFsyLDAsIkIiXSxbMywwLCJCJyJdLFswLDEsIlxcdGF1Il0sWzEsMiwiZyJdLFsyLDMsImYiXV0=
	\begin{tikzcd}
		R & A & B & {B'}
		\arrow["\tau", from=1-1, to=1-2]
		\arrow["g", from=1-2, to=1-3]
		\arrow["f", from=1-3, to=1-4]
	\end{tikzcd} of ring morphisms and if $f$ is an integral morphism and $\ker f$ is a nil ideal then $f^{-1}((A')^*_{B',R'})=A^*_{B,R}$. 
	
	In this section, we use integrality as an important condition so that Lipman diagrams contract. Further, as a consequence of our theorem, we are able to remove the requirement of $\ker f$ to be a nil ideal.
	
	The proof of the next lemma was inspired by Lipman.
	
	\begin{lema}\label{202305241953}
		Let $h: S\rightarrow T$ be an integral morphism of $R$-algebras. 
		
		\begin{enumerate}
			\item [a)] $\ker(h\ten_R h)$ is a nil ideal of $S\ten_RS$;
			
			\item [b)] Suppose that $\ker h$ is a nil ideal of $S$. Then $\overline{I}=h^{-1}(\overline{IT})$, for every $I$ ideal of $S$.
		\end{enumerate}
	\end{lema}
	
	\begin{proof}
		(a) We want to show that $\ker(h\ten_Rh)\sub\bigcap\limits_{\fp\in\tiny{\mbox{Spec}(S\ten_RS)}}\fp$, so it suffices to show that the induced map $(h\ten_Rh)^\sharp:\spec(T\ten_RT)\rightarrow \spec(S\ten_RS)$ is surjective. 
		
		Let $\fp\in\spec(S\ten_RS)$. The domain $\dfrac{S\ten_RS}{\fp}$ can be embedded into an algebraically closed field $F$ such that the kernel of the composition $$S\ten_RS\overset{\mbox{\tiny{projection}}}{\longrightarrow}\dfrac{S\ten_RS}{\fp}\hookrightarrow F$$
		
		\noindent is $\fp$. Let $\alpha$ be this composition and let $\gamma_1,\gamma_2:S\rightarrow S\ten_RS$ be the canonical maps which takes $s\ten_R1\overset{\gamma_1}{\mapsfrom}s\overset{\gamma_2}{\mapsto}1\ten_Rs$, for all $s\in S$. Since $h$ is an integral morphism then there exist ring morphisms $\delta_1,\delta_2:T\rightarrow F$ for which the diagram % https://q.uiver.app/#q=WzAsNixbMSwyLCJTXFx1bmRlcnNldHtSfXtcXG90aW1lc31TIl0sWzMsMiwiRiJdLFswLDEsIlMiXSxbMCwzLCJTIl0sWzAsMCwiVCJdLFswLDQsIlQiXSxbMiw0LCJoIiwyXSxbMiwwLCJcXGdhbW1hXzEiLDJdLFszLDAsIlxcZ2FtbWFfMiJdLFszLDUsImgiLDJdLFswLDEsIlxcYWxwaGEiXSxbNCwxLCJcXGRlbHRhXzEiLDIseyJzdHlsZSI6eyJib2R5Ijp7Im5hbWUiOiJkYXNoZWQifX19XSxbNSwxLCJcXGRlbHRhXzIiLDIseyJzdHlsZSI6eyJib2R5Ijp7Im5hbWUiOiJkYXNoZWQifX19XV0=
		\[\begin{tikzcd}
			T \\
			S \\
			& {S\underset{R}{\otimes}S} && F \\
			S \\
			T
			\arrow["h"', from=2-1, to=1-1]
			\arrow["{\gamma_1}"', from=2-1, to=3-2]
			\arrow["{\gamma_2}", from=4-1, to=3-2]
			\arrow["h"', from=4-1, to=5-1]
			\arrow["\alpha", from=3-2, to=3-4]
			\arrow["{\delta_1}"', dashed, from=1-1, to=3-4]
			\arrow["{\delta_2}"', dashed, from=5-1, to=3-4]
		\end{tikzcd}\]\noindent commutes. The universal property of the tensor product ensures the existence of a unique ring morphism $\beta: T\ten_RT\rightarrow F$ such that $\beta(u\ten_Rv)=\delta_1(u)\delta_2(v)$, $\forall u,v\in T$. Defines $\fq:=\ker\beta\in\spec(T\ten_RT)$. Clearly $\alpha=\beta\circ(h\ten_Rh)$, and since $\fp=\ker\alpha$ then we conclude that $(h\ten_Rh)^{-1}(\fq)=\fp$. 
		
		(b) The persistence of the integral closure of ideals implies $h(\overline{I})\sub\overline{IT}$. Conversely, assume that $x\in h^{-1}(\overline{IT})$. Thus, $y:=h(x)$ is integral over $IT$, and then $yX$ is integral over $T[(IT)X]$, which is integral over $h(S)[h(I)X]$, once $h$ is an integral morphism. Thus, for each there exist $a_i\in I^i$, $i\in\{1,\hdots, n\}$ such that $$h(a_n)+\cdots+h(a_1)y^{n-1}+y^n=0.$$ 
		
		\noindent Hence, $a_n+\cdots+a_1x^{n-1}+x^n\in \ker h\sub \sqrt{(0)}$, which implies the existence of $r\in\bN$ such that $$(a_n+\cdots+a_1x^{n-1}+x^n)^r=0.$$
		
		Therefore, $x\in\overline{I}$.
		
	\end{proof}
	
	Now we state the main theorem of this section.
	
	\begin{teo}[Generalized integral contraction]\label{202306021356}	Suppose that
		\[\begin{tikzcd}
			R & A & B \\
			{R'} & {A'} & {B'}
			\arrow["\tau", from=1-1, to=1-2]
			\arrow["g", from=1-2, to=1-3]
			\arrow["{\tau'}", from=2-1, to=2-2]
			\arrow["{g'}", from=2-2, to=2-3]
			\arrow["{f_R}"', from=1-1, to=2-1]
			\arrow["{f}", from=1-3, to=2-3]
			\arrow["{f_A}", from=1-2, to=2-2]
		\end{tikzcd}\]
		
		\noindent is a Lipman diagram. If $f$ is integral and $\ker p$ is a nil ideal then $$f^{-1}((A')^*_{B',R'})=A^*_{B,R}.$$
	\end{teo}
	
	\begin{proof}
		Since $f$ is integral then $f\ten_Rf$ is integral. By Lemma \ref{202305241953} (a) $\ker(f\ten_Rf)$ is a nil ideal, i.e, $\ker(f\ten_Rf)\sub\sqrt{(0)}$. Further, $p$ is surjective, so it is integral which implies that $\bar{f}:=p\circ(f\ten_Rf)$ is integral. Besides, using that $\ker p$ is a nil ideal, we have $$\ker\bar{f}=(f\ten_Rf)^{-1}(\ker p)\sub(f\ten_Rf)^{-1}(\sqrt{(0)})\sub\sqrt{(f\ten_Rf)^{-1}((0))}$$$$\sub\sqrt{\ker(f\ten_Rf)}\sub\sqrt{\sqrt{(0)}}=\sqrt{(0)}.$$
		
		Thus, $\ker\bar{f}$ is a nil ideal, and Lemma \ref{202305241953} (b) implies that $\bar{f}^{-1}(\overline{\ker\varphi(B'\ten_{R'}B')})=\overline{\ker\varphi}$. From now on, it suffices to proceed exactly as we did in the proof of Theorem \ref{202306011654}.	  
	\end{proof}
	
	As a consequence, we get the Lipman sequence contracts without any requirement about the kernel of $f$.
	
	\begin{cor}[Integral contraction]\label{202305241952}
		Let $R\overset{\tau}{\rightarrow} A \overset{g}{\rightarrow}B\overset{f}{\rightarrow}B'$ be a sequence of ring morphisms, and assume that $f$ is integral.  Then: $$A^*_{B,R}=f^{-1}(A^*_{B',R}).$$
		
	\end{cor}
	
	\begin{proof}
		In Corollary \ref{202306011657} we have seen the sequence gives rise to a Lipman diagram % https://q.uiver.app/#q=WzAsNixbMCwwLCJSIl0sWzEsMCwiQSJdLFsyLDAsIkIiXSxbMCwxLCJSIl0sWzEsMSwiQSJdLFsyLDEsIkInIl0sWzAsMSwiXFx0YXUiXSxbMSwyLCJnIl0sWzMsNCwiXFx0YXUiLDJdLFs0LDUsImZcXGNpcmMgZyIsMl0sWzAsMywiXFx0ZXh0cm17aWR9X1IiLDJdLFsxLDQsIlxcdGV4dHJte2lkfV9BIiwyXSxbMiw1LCJmIiwyXV0=
		\[\begin{tikzcd}
			R & A & B \\
			R & A & {B'}
			\arrow["\tau", from=1-1, to=1-2]
			\arrow["g", from=1-2, to=1-3]
			\arrow["\tau"', from=2-1, to=2-2]
			\arrow["{f\circ g}"', from=2-2, to=2-3]
			\arrow["{\textrm{id}_R}"', from=1-1, to=2-1]
			\arrow["{\textrm{id}_A}"', from=1-2, to=2-2]
			\arrow["f"', from=1-3, to=2-3]
		\end{tikzcd}\].
		
		In this case, $p$ is the identity morphism, and in particular, its kernel is a nil ideal, which finishes the proof once $f$ is integral.
	\end{proof}
	
	\begin{example}\normalfont
		Suppose that $R\overset{\tau}{\rightarrow} A \overset{g}{\rightarrow}B\overset{f=\pi}{\rightarrow}\dfrac{B}{\fm}$ is a sequence of ring morphisms, where $\fm$ is a maximal ideal of $B$, $f=\pi$ is the quotient map and $B$ is not a local ring. Since $\pi$ is surjective then $\pi$ is integral, and Corollary \ref{202305241952} implies that $$\pi^{-1}\left(A^*_{\frac{B}{\fm},R}\right)=A^*_{B,R},$$
		
		\noindent although $\ker f=\fm$ is not a nil ideal of $B$, once if it was, we would conclude that $\spec B=\{\fm\}$, i.e, $B$ would be a local ring.
	\end{example}
	
	As a consequence of the integral contraction, the Lipschitz saturation commutes with quotient by an ideal. This is what the next result says.
	
	\begin{prop}\label{202302061355}
		Let $R\overset{\tau}{\rightarrow} A \overset{g}{\rightarrow}B$ be a sequence of ring morphisms and let $I$ be an ideal of $A$. Then: $$\dfrac{A^*_{B,R}}{IA^*_{B,R}}\cong \left(\dfrac{A}{I}\right)^*_{\frac{B}{IB},R}$$
		
		\noindent through an isomorphism which takes $u+IA^*_{B,R}\mapsto u+IB, \forall u\in A^*_{B,R}$.
	\end{prop}
	
	\begin{proof}
		Consider the diagram % https://q.uiver.app/#q=WzAsNixbMCwwLCJSIl0sWzIsMCwiQSJdLFs0LDAsIkIiXSxbMCwyLCJSIl0sWzIsMiwiXFxkZnJhY3tBfXtJfSJdLFs0LDIsIlxcZGZyYWN7Qn17SUJ9Il0sWzAsMSwiXFx0YXUiXSxbMSwyLCJnIl0sWzMsNCwiXFxwaVxcY2lyY1xcdGF1IiwyXSxbNCw1LCJcXGJhcntnfSIsMl0sWzAsMywiXFx0ZXh0cm17aWR9X1IiLDJdLFsxLDQsIlxccGkiXSxbMiw1LCJcXGJhcntcXHBpfSJdXQ==
		\[\begin{tikzcd}
			R && A && B \\
			\\
			R && {\dfrac{A}{I}} && {\dfrac{B}{IB}}
			\arrow["\tau", from=1-1, to=1-3]
			\arrow["g", from=1-3, to=1-5]
			\arrow["\pi\circ\tau"', from=3-1, to=3-3]
			\arrow["{\bar{g}}"', from=3-3, to=3-5]
			\arrow["{\textrm{id}_R}"', from=1-1, to=3-1]
			\arrow["\pi", from=1-3, to=3-3]
			\arrow["{\bar{\pi}}", from=1-5, to=3-5]
		\end{tikzcd} \eqno{(\star)}\]
		
		\noindent where $\pi$ and $\bar{\pi}$ are the quotient maps, and $\bar{g}$ is the canonical morphism induced by $g$ and $I$. It is clear this diagram is commutative, and since $\pi$ is surjective, Proposition \ref{202306011826} ensures that $(\star)$ is a strong Lipman diagram. In this case, $p:\dfrac{B}{IB}\ten_R\dfrac{B}{IB}\rightarrow \dfrac{B}{IB}\ten_R\dfrac{B}{IB}$ is the identity map, so its kernel is a nil ideal. Besides, the quotient map $\bar{\pi}$ is surjective, so it is integral. By Theorem \ref{202306021356} we have $$\bar{\pi}^{-1}\left(\dfrac{A}{I}\right)^*_{\frac{B}{IB},R}=A^*_{B,R}.$$ 
		
		Again, once $\pi$ is surjective, then the last equation implies that $\bar{\pi}(A^*_{B,R})= \left(\dfrac{A}{I}\right)^*_{\frac{B}{IB},R}$ and the morphism $$\begin{matrix}
			A^*_{B,R} & \longrightarrow & \left(\dfrac{A}{I}\right)^*_{\frac{B}{IB},R}\\
			u & \longmapsto & \bar{\pi}(u)=u+IB
		\end{matrix}$$
		
		\noindent is surjective, whose kernel is $A^*_{B,R}\cap IB=IA^*_{B,R}$. Therefore, this morphism induces the desired isomorphism.
	\end{proof}
	
	Next, we conclude that the quotient inherits to be saturated.
	
	\begin{cor}
		Let $R\overset{\tau}{\rightarrow} A \overset{g}{\rightarrow}B$ be a sequence of ring morphisms and let $I$ be an ideal of $A$. If $A$ is $R$-saturated in $B$ then $\dfrac{A}{I}$ is $R$-saturated in $\dfrac{B}{IB}$.
	\end{cor}
	
	\begin{proof}
		Consider the diagram % https://q.uiver.app/#q=WzAsNixbMCwwLCJSIl0sWzEsMCwiQSJdLFsyLDAsIkIiXSxbMCwxLCJSIl0sWzEsMSwiXFxkZnJhY3tBfXtJfSJdLFsyLDEsIlxcZGZyYWN7Qn17SUJ9Il0sWzAsMSwiXFx0YXUiXSxbMSwyLCJnIl0sWzAsMywiXFx0ZXh0cm17aWR9X1IiLDJdLFsxLDQsIlxccGkiXSxbMiw1LCJcXGJhcntcXHBpfSJdLFszLDQsIlxccGlcXGNpcmNcXHRhdSIsMl0sWzQsNSwiXFxiYXJ7Z30iLDJdXQ==
		\[\begin{tikzcd}
			R & A & B \\
			R & {\dfrac{A}{I}} & {\dfrac{B}{IB}}
			\arrow["\tau", from=1-1, to=1-2]
			\arrow["g", from=1-2, to=1-3]
			\arrow["{\textrm{id}_R}"', from=1-1, to=2-1]
			\arrow["\pi", from=1-2, to=2-2]
			\arrow["{\bar{\pi}}", from=1-3, to=2-3]
			\arrow["\pi\circ\tau"', from=2-1, to=2-2]
			\arrow["{\bar{g}}"', from=2-2, to=2-3]
		\end{tikzcd}\] \noindent and let $\pi^*:\dfrac{A^*_{B,R}}{IA^*_{B,R}}\rightarrow \left(\dfrac{A}{I}\right)^*_{\frac{B}{IB},R}$ be the isomorphism obtained in the previous proposition. 
		
		Let us prove that $\left(\dfrac{A}{I}\right)^*_{\frac{B}{IB},R}\sub\bar{g}\left(\dfrac{A}{I}\right)$. If $w\in \left(\dfrac{A}{I}\right)^*_{\frac{B}{IB},R}$ then we can write $$w=\pi^*(u+IA^*_{B,R})=u+IB,$$ \noindent for some $u\in A^*_{B,R}$. Since $A$ is $R$-saturated in $B$ then $A^*_{B,R}=g(A)$, and then there exists $a\in A$ such that $u=g(a)$. Hence, $w=g(a)+IB=\bar{pi}(g(a))=\bar{g}(\pi(a))=\bar{g}(a+I)\in\bar{g}\left(\dfrac{A}{I}\right)$. 
		
		Therefore, $\dfrac{A}{I}$ is $R$-saturated in $\dfrac{B}{IB}$.

	\end{proof}

	Consider again the sequence of ring morphism $R\overset{\tau}{\rightarrow} A \overset{g}{\rightarrow}B$ where $I$ is an ideal of $A$. We can form another type of diagram % https://q.uiver.app/#q=WzAsNixbMCwwLCJSIl0sWzIsMCwiQSJdLFs0LDAsIkIiXSxbMCwyLCJcXGRmcmFje1J9e1xcdGF1XnstMX0oSSl9Il0sWzIsMiwiXFxkZnJhY3tBfXtJfSJdLFs0LDIsIlxcZGZyYWN7Qn17SUJ9Il0sWzAsMSwiXFx0YXUiXSxbMSwyLCJnIl0sWzMsNCwiXFxiYXJ7XFx0YXV9IiwyXSxbNCw1LCJcXGJhcntnfSIsMl0sWzAsMywiXFxwaV9SIiwyXSxbMSw0LCJcXHBpIl0sWzIsNSwiXFxiYXJ7XFxwaX0iXV0=
	\[\begin{tikzcd}
		R && A && B \\
		\\
		{\dfrac{R}{\tau^{-1}(I)}} && {\dfrac{A}{I}} && {\dfrac{B}{IB}}
		\arrow["\tau", from=1-1, to=1-3]
		\arrow["g", from=1-3, to=1-5]
		\arrow["{\bar{\tau}}"', from=3-1, to=3-3]
		\arrow["{\bar{g}}"', from=3-3, to=3-5]
		\arrow["{\pi_R}"', from=1-1, to=3-1]
		\arrow["\pi", from=1-3, to=3-3]
		\arrow["{\bar{\pi}}", from=1-5, to=3-5]
	\end{tikzcd}\].
	
	Next, we show that we can replace $R$ by $\dfrac{R}{I}$ on the right-hand side of the conclusion of Proposition \ref{202302061355}.
	
	\begin{prop}
		With the above notation, $$\left(\dfrac{A}{I}\right)^*_{\frac{B}{IB},\frac{R}{\tau^{-1}(I)}}=\left(\dfrac{A}{I}\right)^*_{\frac{B}{IB},R}\cong\dfrac{A^*_{B,R}}{IA^*_{B,R}}$$
	\end{prop}
	
	\begin{proof}
		In this case the canonical morphism $\gamma:\dfrac{R}{\tau^{-1}(I)}\ten_R\dfrac{R}{\tau^{-1}(I)}\rightarrow \dfrac{R}{\tau^{-1}(I)}$ is an isomorphism, so its kernel is a nil ideal. Thus\footnote{See Prop. 3.7.1, p. 246 in \cite{EGA}} $\dfrac{R}{\tau^{-1}(I)}$ is a radicial $R$-algebra (with the canonical structure). From the sequence % https://q.uiver.app/#q=WzAsNCxbMCwwLCJSIl0sWzEsMCwiXFxkZnJhY3tSfXtcXHRhdV57LTF9KEkpfSJdLFsyLDAsIlxcZGZyYWN7QX17SX0iXSxbMywwLCJcXGRmcmFje0J9e0lCfSJdLFswLDFdLFsxLDIsIlxcYmFye1xcdGF1fSJdLFsyLDMsIlxcYmFye2d9Il1d
		\[\begin{tikzcd}
			R & {\dfrac{R}{\tau^{-1}(I)}} & {\dfrac{A}{I}} & {\dfrac{B}{IB}}
			\arrow[from=1-1, to=1-2]
			\arrow["{\bar{\tau}}", from=1-2, to=1-3]
			\arrow["{\bar{g}}", from=1-3, to=1-4]
		\end{tikzcd}\]
		
		\noindent Proposition 1.3 of \cite{lipman1} ensures that $\left(\dfrac{A}{I}\right)^*_{\frac{B}{IB},\frac{R}{\tau^{-1}(I)}}=\left(\dfrac{A}{I}\right)^*_{\frac{B}{IB},R}$.
	\end{proof}

Now we apply the main theorem of this section to get a new class of algebras where the conclusion of Proposition 1.3 of \cite{lipman1} is satisfied.

\begin{prop}
	Consider a sequence of ring morphisms % https://q.uiver.app/#q=WzAsNCxbMCwwLCJSIl0sWzEsMCwiUiciXSxbMiwwLCJBIl0sWzMsMCwiQiJdLFswLDEsImZfUiJdLFsxLDJdLFsyLDMsImciXV0=
	$\begin{tikzcd}
		R & {R'} & A & B
		\arrow[from=1-1, to=1-2]
		\arrow[from=1-2, to=1-3]
		\arrow["g", from=1-3, to=1-4]
	\end{tikzcd}$. If $\ker(B\ten_RB\rightarrow B\ten_{R'}B)$ (which is $\ker p$) is a nil ideal then $$A^*_{B,R}=A^*_{B,R'}.$$ 
\end{prop}

\begin{proof}
	From the above sequence we get a strong Lipman diagram % https://q.uiver.app/#q=WzAsNixbMCwwLCJSIl0sWzIsMCwiQiJdLFswLDEsIlInIl0sWzEsMSwiQSJdLFsyLDEsIkIiXSxbMSwwLCJBIl0sWzAsNV0sWzAsMiwiZl9SIiwyXSxbMiwzXSxbMyw0LCJnIiwyXSxbNSwxLCJnIl0sWzUsMywiXFx0ZXh0cm17aWR9X0EiXSxbMSw0LCJcXHRleHRybXtpZH1fQiJdXQ==
	\[\begin{tikzcd}
		R & A & B \\
		{R'} & A & B
		\arrow[from=1-1, to=1-2]
		\arrow["{f_R}"', from=1-1, to=2-1]
		\arrow[from=2-1, to=2-2]
		\arrow["g"', from=2-2, to=2-3]
		\arrow["g", from=1-2, to=1-3]
		\arrow["{\textrm{id}_A}", from=1-2, to=2-2]
		\arrow["{\textrm{id}_B}", from=1-3, to=2-3]
	\end{tikzcd}\]\noindent  and since $\id_B$ is integral and $\ker p$ is a nil ideal then Theorem \ref{202306021356} implies  $\id_B^{-1}(A^*_{B,R'})=A^*_{B,R}$. Hence, $A^*_{B,R}=A^*_{B,R'}$.
\end{proof}

In the following result, we combine integrality and universal injectivity to derive a contraction result.

\begin{prop}
Suppose that
\[\begin{tikzcd}
	R & A & B \\
	{R'} & {A'} & {B'}
	\arrow["\tau", from=1-1, to=1-2]
	\arrow["g", from=1-2, to=1-3]
	\arrow["{\tau'}", from=2-1, to=2-2]
	\arrow["{g'}", from=2-2, to=2-3]
	\arrow["{f_R}"', from=1-1, to=2-1]
	\arrow["{f}", from=1-3, to=2-3]
	\arrow["{f_A}", from=1-2, to=2-2]
\end{tikzcd}\]

\noindent is a Lipman diagram. If $f_R$ is universally injective and $f$ is integral then $$f^{-1}((A')^*_{B',R'})=A^*_{B,R}.$$
\end{prop}
	
\begin{proof}
	Since $f_R$ is universally injective then Lemma \ref{202306082225} implies that $p$ is an isomorphism. In particular, $\ker p=(0)$ is a nil ideal. From now on, the result is a direct consequence of Theorem \ref{202306021356}.
\end{proof}

	\bibliographystyle{plain}
	\bibliography{referencias}

\end{document}